\documentclass[letterpaper, 10 pt, conference]{ieeeconf}  

\IEEEoverridecommandlockouts                              
\usepackage{amsmath,amsfonts,amssymb} 
\usepackage{tikz}
\usepackage{algorithm}
\usepackage{algorithmic}
\usepackage{hyperref}
\newtheorem{definition}{Definition}
\newcommand\norm[1]{\lVert#1\rVert}
\newcommand{\minimize}[1]{\underset{{#1}}{\text{minimize}}}

\newcommand{\st}{\text{subject to}}
\newcommand{\mb}[1]{\mathbf{#1}}
\newcommand{\mbg}[1]{\boldsymbol{#1}}
\newcommand{\theoremsymbol}{\hfill \ensuremath{\triangleleft}}
\def\Plus{\texttt{+}}

\title{\LARGE \bf
Optimization over Trained Neural Networks: Difference-of-Convex Algorithm and Application to Data Center Scheduling
}

\author{Xinwei Liu and Vladimir Dvorkin%
\thanks{Xinwei Liu and Vladimir Dvorkin are with the Department of Electrical Engineering and Computer Science, University of Michigan, MI 48109, USA. E-mail: {\tt \{xinwei,dvorkin\}@umich.edu}}
}

\begin{document}
\begingroup
\allowdisplaybreaks

\maketitle
\thispagestyle{empty}
\pagestyle{empty}

\begin{abstract}
When solving decision-making problems with mathematical optimization, some constraints or objectives may lack analytic expressions but can be approximated from data. When an approximation is made by neural networks, the underlying problem becomes optimization over trained neural networks. Despite recent improvements with cutting planes, relaxations, and heuristics, the problem remains difficult to solve in practice. We propose a new solution based on a bilinear problem reformulation that penalizes ReLU constraints in the objective function. This reformulation makes the problem amenable to efficient difference-of-convex algorithms (DCA), for which we propose a principled approach to penalty selection that facilitates convergence to stationary points of the original problem. We apply the DCA to the problem of the least-cost allocation of data center electricity demand in a power grid, reporting significant savings in congested cases.
\end{abstract}

\section{Introduction}

Optimization over trained neural networks has gained momentum thanks to their remarkable expressivity in representing complex objectives and constraints in decision-making. Rectified Linear Unit (ReLU) neural networks can be trained to learn complex laws from data and then integrated into mathematical optimization as mixed-integer constraints \cite{tjeng2017evaluating,grimstad2019relu}. For example, ReLU networks efficiently approximate natural gas and power flows from data, enabling optimal dispatch of gas systems and power grids \cite{dvorkin2023emission, kody2022modeling}---the application domain of interest to this work. 

For practically sized neural networks, such optimization becomes computationally expensive. The standard mixed-integer programming (MIP) approach enjoys Big-M reformulation of ReLU constraints but suffers from numerical issues when M is large \cite{tjeng2017evaluating,grimstad2019relu}. The SOS1 approach resolves numerical issues \cite{turner2024pyscipoptmlembeddingtrainedmachine} but requires computationally expensive branch-and-cut procedures. At the rate of one binary variable (or SOS1 constraint) per ReLU, such methods quickly become intractable. To improve tractability, the literature suggests identifying stable neurons \cite{xiao2018training} and patterns in neural activations \cite{serra2020empirical}, as well as pruning \cite{serra2020lossless}. Considerable progress has been made in the space of neural network verification, offering improvements via semidefinite programming relaxation of ReLU neurons \cite{dathathri2020enabling,fazlyab2020safety} and bound propagation \cite{wang2021beta}, although their applications outside the verification space have yet to be explored. In addition, a principled heuristic proposed in \cite{tong2024optimization} exploring global and local linear relaxations of neural networks has shown to outperform MIP solvers.

Our interest in optimization over trained neural networks is motivated by the growing electricity demand of AI. Large networks of data centers can help grid operators reduce the cost of grid operations by spatially reallocating electricity demand among electrically remote locations. This problem can be solved using different types of coordination (e.g., see \cite{dvorkin2024agent} for a brief overview), but more naturally via strategic market participation of data centers. By seeking the least-cost demand allocation in the grid, they implicitly help the system achieve the minimum operating cost. Strategic participation models either feature restrictive assumptions (e.g., linear price model \cite{zhao2022strategic}) or require grid data \cite{dvorkin2018consensus} which is not available to regular participants. The use of neural networks provides a more practical solution: They are trained on historical market outcomes to approximate the price model, which can then be used to optimize demand allocations.

\textit{Contributions:} Inspired by their potential impact on grid operations, we develop a new efficient algorithm for optimization over trained neural networks. Our algorithm relaxes the ReLU constraints by penalizing them in the objective function. Such a formulation makes the problem amenable to difference-of-convex programming \cite{an2005dc}, which iteratively solves simple convex sub-problems, thereby avoiding the computational complexity of the standard MIP formulations. However, the challenge is to compute the penalty parameter, which sufficiently penalizes the ReLU constraints without impacting optimality. We establish the principled approach to computing such a penalty. Our approach is inspired by solution techniques for complementarity problems whose structure is similar to the problem at hand \cite{jara2018study}. We provide an application to strategic participation of data centers in electricity markets, where trained neural networks are used to approximate the market price as a function of data center demands. Our simulations on small and large neural networks reveal significant cost savings for data center owners. 

\textit{Paper Organization:} Sec. \ref{sec:preliminaries} provides preliminaries on the problem at hand. Sec. \ref{sec:DCA} introduces the main algorithm and our approach to selecting the penalty. Sec. \ref{sec:application}  and \ref{sec:experiments} provide application and numerical results. Sec. \ref{sec:conclusion} concludes. 

\section{Optimization over Trained Neural Networks}\label{sec:preliminaries}

Consider a \textit{trained} neural network (NN) with the following mathematical representation:
\begin{align*}
    \mbg{\lambda} = \mb{W}^{N+1}\text{ReLU}(\dots\text{ReLU}(\mb{W}^1\mb{d}+\mb{b}^{1})\dots)+\mb{b}^{N+1},
\end{align*}
i.e.,  feedforward network with fully-connected layers of ReLU neurons. Consider the input $\mb{d}=[d_{1},\dots,d_{n}]^{\top}$ from a bounded domain $\mathcal{D}$ and the output $\mbg{\lambda}=[\lambda_{1},\dots,\lambda_{m}]^{\top}$. Each hidden layer $i=1,\dots,N$ has an output $\mb{y}^{i}=[y_{1}^{i},\dots,y_{p}^{i}]^{\top}.$ Let $\mb{W}^{i}$ be the $n_{i}\times n_{i-1}$ matrix of weights connecting layer $i$ and $i-1,$ and $\mb{b}^{i}$ be the bias vector. The output $\mb{y}^{i}$ is given by the ReLU activation function $\text{ReLU}(\mb{W}^{i}\mb{y}^{i-1} + \mb{b}^{i})=\text{max}\{\mb{0},\mb{W}^{i}\mb{y}^{i-1} + \mb{b}^{i}\},$ understood element-wise.      

In optimization over trained neutral networks, the parameters $\mb{W}^{i}$ and $\mb{b}^{i}$ for each layer $i$ are fixed, and the decision variables include the input $\mb{d}$, and the outputs $\mb{y}^{1},\dots,\mb{y}^{N}$ and $\mbg{\lambda}$ of the hidden and final linear layers, respectively. The optimization problem takes the form: 
\begin{subequations}\label{prob:nn_opt_intractable}
\begin{align}
   \minimize{\mb{d} \in \mathcal{D},\mb{y}, \mbg{\lambda}}
   \;\;\!&  c(\mbg{\lambda})\\
    \st\;\;\!
    & \mb{y}^1=\text{ReLU} (\mb{W}^1 \mb{d}+\mb{b}^1),\label{prob:nn_opt_intractable_relu1}\\
    & \mb{y}^i=\text{ReLU} (\mb{W}^i \mb{y}^{i-1}+\mb{b}^i), \forall i=2,\dots,N \label{prob:nn_opt_intractable_relu2}\\
     & \mbg{\lambda}= \mb{W}^{N + 1} \mb{y}^{N}+\mb{b}^{N+1},
    \end{align}
\end{subequations}
where the input $\mb{d}$ belongs to its bounded domain $\mathcal{D}$, $c:\mathbb{R}^{m}\mapsto\mathbb{R}$ is the convex cost function acting on NN'
s outputs, and variable $\mb{y}=[\mb{y}^{1\top},\dots,\mb{y}^{N\top}]^{\top}$ stacks the outputs of hidden layers. The goal of this optimization is to find the optimal NN input $\mb{d}^{\star}$ that minimizes the cost function.

We relate the optimization over trained NNs to mixed-complementarity problems. The neuron activations $y=\text{ReLU}(a)$ are represented by complementarity constraints:
\begin{subequations}
\begin{align}
&y= a + v,\\
 & 0 \leqslant y \perp v \geqslant 0,\label{eq:sos1-2}
\end{align}
\end{subequations}
where $v$ is an auxiliary variable, and expression $y \perp v$ means that $y$ and $v$ are orthogonal, i.e., their inner product is zero. We can verify this ReLU formulation by examining two scenarios, $a \geqslant 0$ and $a \leqslant 0$. If $a \geqslant 0$, then $y = a$, and $v = 0$ and all conditions are satisfied. If $a \leqslant 0$, then $y = 0$, $v = -a \geqslant 0$, again, all conditions are satisfied. We can also regard the last condition \eqref{eq:sos1-2} as a SOS1 constraint \cite{turner2024pyscipoptmlembeddingtrainedmachine}, or address it through Big-M reformulation \cite{tjeng2017evaluating}; in this paper, however, we will take a different path to address this constraint. We rewrite problem \eqref{prob:nn_opt_intractable} as follows:
\begin{subequations}\label{prob:nn_opt_mixed_comp}
\begin{align}
   \minimize{\mb{d} \in \mathcal{D},\mb{y}, \mb{v}, \mbg{\lambda}}
   \;\;\!&  c(\mbg{\lambda})\\
    \st\;\;\!
    & \mb{y}^1=\mb{W}^1 \mb{d}+\mb{b}^1 + \mb{v}^1,\label{prob:nn_opt_intractable_relu1}\\
    & \mb{y}^i=\mb{W}^i \mb{y}^{i-1}+\mb{b}^i + \mb{v}^i, \forall i=2,\dots,N \label{prob:nn_opt_intractable_relu2}\\
    &\mb{0}\leqslant \mb{y}^{i} \perp \mb{v}^{i} \geqslant \mb{0},\quad\quad\quad\!\;\forall i=1,\dots,N \label{prob:nn_opt_comp}\\
     & \mbg{\lambda}= \mb{W}^{N + 1} \mb{y}^{N}+\mb{b}^{N+1},
    \end{align}
\end{subequations}

The number of complementarity constraints is equal to the number of hidden neurons. At the rate of one constraint per ReLU, the problem quickly becomes intractable. In the next section, we introduce an algorithm to address this challenge. 

\section{Difference-of-Convex Algorithm for Optimization Over Trained Neural Networks}\label{sec:DCA}
In this section, we first introduce a bilinear reformulation of problem \eqref{prob:nn_opt_mixed_comp} that penalizes complementarity constraint violations in the objective function. The difference--of--convex algorithm (DCA) is then introduced to solve the bilinear formulation. Leveraging the specifics of trained NNs, we then present a principled method to optimize the DCA parameters to converge to the solution of the original problem.

\subsection{Bilinear Formulation and DCA Algorithm}

The bilinear formulation relaxes the complementarity constraints \eqref{prob:nn_opt_comp} and penalizes them in the objective function:
\begin{subequations}\label{prob:nn_opt_bileinear}
\begin{align}
   \minimize{(\mb{d},\mb{y},\mb{v},\mbg{\lambda})\in\mathcal{O}}
   \;\;\!&  c(\mbg{\lambda}) + \rho\sum_{i=1}^{N}\mb{y}^{i\top}\mb{v}^{i}\label{prob:nn_opt_bileinear_obj}\\
    \st\;\;\!
    & \mb{y}^1=\mb{W}^1 \mb{d}+\mb{b}^1 + \mb{v}^1,\label{prob:nn_opt_bileinear_relu1}\\
    & \mb{y}^i=\mb{W}^i \mb{y}^{i-1}+\mb{b}^i + \mb{v}^i, \;\forall i=2,\dots,N \label{prob:nn_opt_bileinear_relu2}\\
    & \mbg{\lambda}= \mb{W}^{N + 1} \mb{y}^{N}+\mb{b}^{N+1},\label{prob:nn_opt_bileinear_last}
    \end{align}
\end{subequations}
where $\rho>0$ is a penalty parameter and set $\mathcal{O}$ contains the input domain $\mathcal{D}$ and non-negativity of variables $\mb{y}$ and $\mb{v}$. Under such constraints, the penalty function $\rho\mb{y}^{i\top}\mb{v}^{i}$ for each layer $i$ is zero if and only if $\mb{y}^{i} \perp \mb{v}^{i}$, as required by \eqref{prob:nn_opt_comp}. By solving this penalized problem, we aim to obtain the solution to the original non-convex problem \eqref{prob:nn_opt_mixed_comp}. Finding the global optimum solution of \eqref{prob:nn_opt_mixed_comp} is NP-hard, so we seek \textit{strongly stationary solutions}--to be formally defined later--which satisfy the Karush-Kuhn-Tucker conditions (KKTs)  of \eqref{prob:nn_opt_mixed_comp}. We observe that the penalty functions for each layer $i$ can equivalently be written as \cite{jara2018study}: 
\begin{align}\label{eq:bilinear_penalty_eq}
    \rho\mb{y}^{i\top}\mb{v}^{i} = \frac{\rho}{4} \norm{\mb{y}^{i} + \mb{v}^{i}}_2^2 - \frac{\rho}{4}\norm{\mb{y}^{i} - \mb{v}^{i}}_2^2
\end{align}
Substituting \eqref{eq:bilinear_penalty_eq} into the objective function \eqref{prob:nn_opt_bileinear_obj}, we represent the latter as the difference--of--convex function 
\begin{align}
   f(\mb{y},\mb{v})=\underbrace{c(\mbg{\lambda}) + \frac{\rho}{4}\sum_{i=1}^{N}\norm{\mb{y}^{i} + \mb{v}^{i}}_2^2}_{f_1(\mb{y},\mb{v})} 
   - \underbrace{\frac{\rho}{4}\sum_{i=1}^{N}\norm{\mb{y}^{i} - \mb{v}^{i}}_2^2}_{f_2(\mb{y},\mb{v})}\label{eq:dcf}\nonumber\\[-0.5cm]
\end{align}
where the non-convex objective function $f$ is given by the difference of two convex functions $f_1$ and $f_2$. We do not include variable $\mbg{\lambda}$ into the arguments of $f$ and $f_1$, as it can be represented by the affine transformation of the output of the last hidden layer $\mb{y}^{N}$. The basic DCA minimizes \eqref{eq:dcf} over iterations subject to constraints \eqref{prob:nn_opt_bileinear_relu1}--\eqref{prob:nn_opt_bileinear_last}.
It overestimates the concave part of \eqref{eq:dcf} by a linear function using a gradient evaluated at a current iterate and solves a resultant auxiliary (convex) sub-problem to obtain a new iterate \cite{an2005dc}. In terms of function \eqref{eq:dcf}, after computing the gradient of $f_2$, the corresponding sub-problem at iteration $k$ takes the form 
\begin{subequations}\label{prob:bileinear_sub_problem}
\begin{align}
   \minimize{(\mb{d},\mb{y},\mb{v},\mbg{\lambda})\in\mathcal{O}}
   \;\;&  c(\mbg{\lambda}) + \frac{\rho}{4}\sum_{i=1}^{N}\norm{\mb{y}^{i} + \mb{v}^{i}}_2^2 \nonumber\\
   &\quad\quad- \frac{\rho}{2}\sum_{i=1}^{N}(\mb{y}^{i(k)} - \mb{v}^{i(k)})^{\top}(\mb{y}^{i} - \mb{v}^{i})\\
    \st\;\;
    & \mb{y}^1=\mb{W}^1 \mb{d}+\mb{b}^1 + \mb{v}^1,\\
    & \mb{y}^i=\mb{W}^i \mb{y}^{i-1}+\mb{b}^i + \mb{v}^i, \;\forall i=2,\dots,N\\
    & \mbg{\lambda}= \mb{W}^{N + 1} \mb{y}^{N}+\mb{b}^{N+1},
    \end{align}
\end{subequations}
where the variables $\mb{y}^{i(k)},\mb{v}^{i(k)}$ are fixed to the corresponding values from iteration $k$, and the last term represents the gradients of $f_{2}$ evaluated at $(\mb{y}^{i(k)},\mb{v}^{i(k)})$. This problem is convex and easy to solve. The DCA to solve problem \eqref{prob:nn_opt_bileinear} is summarized in Alg. \ref{alg:bilinear}. It requires a feasible initial guess for set $\mathcal{O}$ and constraints \eqref{prob:nn_opt_bileinear_relu1}--\eqref{prob:nn_opt_bileinear_last}, and a tolerance parameter $\varepsilon_{\text{tol}}$. The algorithm outputs the optimized $\mb{d}^{\star}$ that minimizes the objective function of the penalized problem \eqref{prob:nn_opt_bileinear}. 
\begin{algorithm}[t]
\caption{DCA for solving bilinear problem \eqref{prob:nn_opt_bileinear}}
\label{alg:bilinear}
\begin{algorithmic}[1]
\STATE \textbf{input:} feasible guess $\mb{d}^0, \mb{y}^0, \mb{v}^0, \mbg{\lambda}^0$, tolerance $\varepsilon_{\text{tol}}>0$
\STATE \textbf{output:} optimized NN input $\mb{d}^{\star}$
\vspace{5pt}
    \REPEAT
        \STATE set $k \leftarrow k + 1$.
        \STATE get $\mb{d}^{k+1}, \mb{y}^{k+1}, \mb{v}^{k+1}, \mbg{\lambda}^{k+1}$ by solving \eqref{prob:bileinear_sub_problem}
    \UNTIL{$f(\mb{y}^k, \mb{v}^k) - f(\mb{y}^{k+1},\mb{v}^{k+1}) \leq \varepsilon_{\text{tol}}$}
\vspace{5pt}
\STATE \textbf{return} $\mb{d}^{\star} \leftarrow \mb{d}^{k+1}$
\end{algorithmic}
\end{algorithm}

\subsection{Selecting Penalty Parameter $\rho$}\label{subsec:rho}

The penalty parameter $\rho$ in the bilinear problem \eqref{prob:nn_opt_bileinear} plays an important role in retrieving the solution to the original problem \eqref{prob:nn_opt_mixed_comp}. A small $\rho$ does not sufficiently penalize the relaxed complementarity constraints in \eqref{prob:nn_opt_bileinear_obj}, while a large penalty prioritizes complementarity constraint satisfaction over cost minimization. Moreover, a large penalty slows down DCA convergence, as we illustrate later in Sec. \ref{sec:application}. 

Hence, we seek the smallest $\rho$ that satisfies the complementarity constraints without impacting convergence. Our analysis follows that underlying \cite[Prop. 6]{jara2018study} while relying on specifics of trained NNs. The optimal $\rho$ highly depends on the NN at hand, so we pose the following important question: given a trained NN, what is the penalty $\rho$ that drives the DCA to a stationary solution of problem \eqref{prob:nn_opt_mixed_comp}?

We first rewrite problem \eqref{prob:nn_opt_mixed_comp} in a more general form:
\begin{subequations}\label{prob:general}
\begin{align}
   \underset{\mb{d},\mb{y},\mb{v}}{\text{minimize}}\quad& \mb{c}^\top \mb{y}
   \label{general_obj}\\
    \st\quad
    & \mb{A}\mb{d} \geqslant \mb{f}, \label{general_con_1}\\
    & \mb{V}\mb{d} + \mb{W}\mb{y} + \mb{b} = \mb{v}, \label{general_con_2}\\
    &\mb{0}\leqslant \mb{y}^{i} \perp \mb{v}^{i} \geqslant \mb{0},\quad\quad\quad\!\;\forall i=1,\dots,N,\label{general_con_3}
\end{align}
\end{subequations}
where the objective represents the output of the neural network we optimize for. The NN input $\mb{d}$  is constrained by inequality \eqref{general_con_1}, and equality \eqref{general_con_2} represents NN constraints \eqref{prob:nn_opt_bileinear_relu1}--\eqref{prob:nn_opt_bileinear_relu2}. The bilinear formulation then becomes
\begin{subequations}\label{eq:bilinear_general}
\begin{align}
   \underset{\mb{d},\mb{y},\mb{v}}{\text{minimize}}\quad&\mb{c}^\top \mb{y} + \rho \mb{y}^\top \mb{v}\\
    \st\quad
    & \mb{A}\mb{d} \geqslant \mb{f}, \\
    & \mb{V}\mb{d} + \mb{W}\mb{y} + \mb{b} = \mb{v}, \\
    & \mb{y}_i, \mb{v}_i \geqslant 0, \quad\quad\quad\!\;\forall i=1,\dots,N 
\end{align}
\end{subequations}

To compute $\rho$, we explore the relationship between several optimization problems, as shown in Fig. \ref{fig:relationships}. First, we compute a stationary point of problem \eqref{prob:general}---the point which satisfies its KKTs---by solving a relaxed version of \eqref{prob:general}, whose solution is still feasible for \eqref{prob:general}.  Then, we find $\rho$ by relating the KKTs of the relaxation of \eqref{prob:general} and KKTs of the bilinear formulation \eqref{eq:bilinear_general}, both enforced on the stationary point. 
\begin{figure}
    \centering\small
    \begin{tikzpicture}
    \node[align=center,draw=black, rounded corners = 5pt,fill=gray!10] (original) at (0,0) {Strongly stationary point of problem \eqref{prob:general}};

    \node[align=center,draw=black, rounded corners = 5pt,fill=gray!10] (RxLP) at (0,-1.05) {Solution to relaxation \eqref{prob:general_relax} with empty set $\mathcal{I}_0$ };


    \node[align=center,draw=black, rounded corners = 5pt,fill=gray!10] (Bilinear_2) at (0,-2.25) {DCA solution to bilinear problem \eqref{eq:bilinear_general} for $\rho>\overline{\rho}$ \\ if complementarity $\mb{0}\leqslant\overline{\mb{y}}\perp\overline{\mb{v}}\geqslant\mb{0}$ is satisfied};

    \draw[double,<->,>=stealth,black,line width = 0.025cm] (original) -- (RxLP);


    \draw[double,<->,>=stealth,black,line width = 0.025cm] (RxLP) -- (Bilinear_2);
\end{tikzpicture}
    \caption{Relationships between stationary points of problems in Sec. \ref{subsec:rho}. Such points satisfy the same set of KKT conditions \cite[Prop.6]{jara2018study}.}
    \label{fig:relationships}
\end{figure}
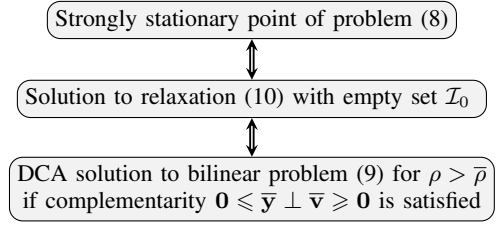

For the review of stationarity concepts, we refer to \cite{fletcher2006local,jara2018study}; here, we are interested in the strong stationarity. 

\vspace{0.1cm}
\begin{definition}[Strongly stationary point]
Let $(\tilde{\mb{d}}, \tilde{\mb{y}},\tilde{\mb{v}})$ be some feasible point for constraints \eqref{general_con_1}--\eqref{general_con_3} and consider the following problem as a relaxation of problem \eqref{prob:general}:
\begin{subequations}\label{prob:general_relax}
\begin{align}
   \underset{\mb{d},\mb{y},\mb{v}}{\text{minimize}}\quad&\mb{c}^\top \mb{y}\\
    \st\quad
    & \mb{A}\mb{d} \geqslant \mb{f}\quad \colon\mbg{\lambda}, \label{general_relax_con_1}\\
    & \mb{V}\mb{d} + \mb{W}\mb{y} + \mb{b} = \mb{v}\quad \colon\mbg{\mu}, \label{general_relax_con_2}\\
    & y_i = 0\quad\colon\mu_i^{y_{0}}  \quad \forall i \in \mathcal{I}_y (\tilde{\mb{y}},\tilde{\mb{v}}),\label{general_relax_con_3}\\
    & y_i \geqslant 0\quad\colon\mu_i^{y_{\Plus}} \quad \forall i \in \mathcal{I}_v (\tilde{\mb{y}},\tilde{\mb{v}}),\label{general_relax_con_4}\\
    & v_i = 0\quad\colon\mu_i^{v_{0}} \quad \forall i \in \mathcal{I}_v (\tilde{\mb{y}},\tilde{\mb{v}}),\label{general_relax_con_5}\\
    & v_i \geqslant 0\quad\colon\mu_i^{v_{\Plus}} \quad \forall i \in \mathcal{I}_y (\tilde{\mb{y}},\tilde{\mb{v}}),\label{general_relax_con_6}\\
    & y_i, v_i \geqslant 0, \quad \quad\;\;\!  \forall i \in \mathcal{I}_0(\tilde{\mb{y}},\tilde{\mb{v}}),\label{general_relax_con_7}
\end{align}
\end{subequations}
whose dual variables are stated after the colon signs, and the index sets are defined as follows
\begin{align*}
\mathcal{I}_y(\tilde{\mb{y}},\tilde{\mb{v}}) &\triangleq\left\{i \mid \tilde{y}_i=0<\tilde{v}_i\right\},\\
\mathcal{I}_w(\tilde{\mb{y}},\tilde{\mb{v}}) &\triangleq\left\{i \mid \tilde{y}_i>0=\tilde{v}_i\right\}, \\
\mathcal{I}_0(\tilde{\mb{y}},\tilde{\mb{v}}) &\triangleq\left\{i \mid \tilde{y}_i=0=\tilde{v}_i\right\}.
\end{align*}
The solution $(\overline{\mb{d}}, \overline{\mb{y}},\overline{\mb{v}})$ of the relaxed problem \eqref{prob:general_relax} is called a \textit{strongly stationary point} of \eqref{prob:general} if it is feasible for \eqref{prob:general}. \theoremsymbol
\end{definition}
\vspace{0.1cm}


To compute a strongly stationary point, we sample feasible points $(\tilde{\mb{d}}, \tilde{\mb{y}},\tilde{\mb{v}})$ in expectation that the solution of the relaxed problem \eqref{prob:general_relax} is also feasible for \eqref{prob:general}. Conveniently, for trained NNs, the feasible point $(\tilde{\mb{d}}, \tilde{\mb{y}},\tilde{\mb{v}})$ can be obtained by sampling a random $\tilde{\mb{d}}$ from the input domain $\mathcal{D}$ and propagating it through the network, taking the note of corresponding variables $\tilde{\mb{y}}$ and $\tilde{\mb{v}}$ for each neuron. For active neurons, $\tilde{y}_i$ takes the neuron output and $\tilde{v}_i=0$. For inactive neurons, $\tilde{y}_i=0$ and $\tilde{v}_i$ takes the negative of the neuron's input. To ensure that the relaxed solution is feasible for \eqref{prob:general}, we look for such point $(\tilde{\mb{d}}, \tilde{\mb{y}},\tilde{\mb{v}})$ which renders index set $\mathcal{I}_0(\tilde{\mb{y}},\tilde{\mb{v}})$ empty, or equivalently $\tilde{y}_i \neq \tilde{v}_i, \forall i=1,\dots,N$. Neurons that fall into this index set feature zero input and zero output, which is rare, making the search very fast. Moreover, if there are many neurons with zero outputs, it is a good indication that the network can be losslessly reduced \cite{serra2020lossless}. Once such point is found, we solve relaxation \eqref{prob:general_relax} on $(\tilde{\mb{d}}, \tilde{\mb{y}},\tilde{\mb{v}})$ to get the strongly stationary point $(\overline{\mb{d}}, \overline{\mb{y}},\overline{\mb{v}})$. 

On the one hand, the subset of KKTs of \eqref{prob:general_relax} include
\begin{subequations}\label{KKTs_relax}
\begin{align}
& \mb{A}^\top \mbg{\lambda}+\mb{V}^\top \mbg{\mu}^v=\mb{0}, \label{KKTs_relax_1}\\
& \mb{W}^\top \mbg{\mu}^v+\mbg{\mu}^y=\mb{c}, \label{KKTs_relax_2}\\
& \mb{0} \leqslant \mbg{\lambda} \perp  \mb{A}\mb{d} - \mb{f} \geqslant \mb{0}, \label{KKTs_relax_3}\\
& \mb{y}^\top \mbg{\mu}^y = 0,
\;\;\;  \mb{v}^\top \mbg{\mu}^v = 0, \label{KKTs_relax_4}\\
& \mu^y_i  \geqslant 0, \quad \forall i \in \mathcal{I}_v(\tilde{\mb{y}},\tilde{\mb{v}}) \label{KKTs_relax_5}\\
& \mu^v_i \geqslant 0, \quad \forall i \in \mathcal{I}_y(\tilde{\mb{y}},\tilde{\mb{v}}) \label{KKTs_relax_6}
\end{align}
\end{subequations}
where, $\mbg{\mu}^v=[\mbg{\mu}^{v_{0}\top}\;\!\mbg{\mu}^{v_{\Plus}\top}]^{\top}$ and $\mbg{\mu}^y=[\mbg{\mu}^{y_{0}\top}\;\!\mbg{\mu}^{y_{\Plus}\top}]^{\top}$, and constraint \eqref{general_relax_con_7} is disregarded as set $\mathcal{I}_0$ is empty. On the other hand, the KKTs of the bilinear problem \eqref{eq:bilinear_general} include
\begin{subequations}\label{KKTs_bilinear}
\begin{align}
\quad\quad\quad\quad\quad\;\;& \mb{A}^\top \mbg{\lambda}+\mb{V}^\top \mbg{\mu}=\mb{0}, \label{KKTs_bilinear_1}\\
& \mb{0} \leqslant \mb{y} \perp-\mb{W}^\top \mbg{\mu}+ \mb{c} + \rho\mb{v} \geqslant \mb{0},\label{KKTs_bilinear_2}\\
& \mb{0} \leqslant\mbg{\lambda} \perp  \mb{A}\mb{d} - \mb{f} \geqslant \mb{0},\label{KKTs_bilinear_3}\\
& \mb{0} \leqslant \mb{v} \perp \rho \mb{y} + \mbg{\mu} \geqslant \mb{0}.\label{KKTs_bilinear_4}
\end{align}
\end{subequations}

Now we pose the following question: for which $\rho$, does the KKT point of \eqref{KKTs_relax} satisfy the KKTs in \eqref{KKTs_bilinear}? Consider a strongly stationary point $(\overline{\mb{d}}, \overline{\mb{y}},\overline{\mb{v}})$, and the corresponding triplet of dual variables $(\overline{\mbg{\lambda}}, \overline{\mbg{\mu}}^v,\overline{\mbg{\mu}}^y),$ jointly satisfying KKTs \eqref{KKTs_relax}. Substituting $\mb{c}$ from \eqref{KKTs_relax_2} to \eqref{KKTs_bilinear_2}, yields the first bound $\rho\geqslant-\overline{\mbg{\mu}}^y/\overline{\mb{v}}$  at the stationary point. Similarly, from the \eqref{KKTs_bilinear_4}, we have the second bound $\rho\geqslant-\overline{\mbg{\mu}}^v/\overline{\mb{y}}$. The two bounds are respected simultaneously for $\rho\geqslant\overline{\rho}$, where 
\begin{align}
    \bar{\rho} \triangleq \max \left\{0,\left\{\left.-\frac{\overline{\mu}_i^y}{\bar{v}_i} \right\rvert\, \bar{v}_i>0\right\},\left\{\left.-\frac{\overline{\mu}_i^v}{\bar{y}_i} \right\rvert\, \bar{y}_i>0\right\}\right\}.\label{eq:rho_bar}
\end{align}
That is, the strongly stationary point $(\overline{\mb{d}},\overline{\mb{y}},\overline{\mb{v}})$ is stationary for the bilinear problem \eqref{eq:bilinear_general} when $\rho \geqslant \bar{\rho}$. Unfortunately, it does not guarantee the DCA convergence to the point satisfying the complementarity condition $\mb{0}\leqslant\overline{\mb{y}}\perp\overline{\mb{v}}\geqslant\mb{0}$. However, if the DCA solution satisfies this condition, it returns a strongly stationary point of problem \eqref{prob:general}.

The overall search logic is summarized in Alg. \ref{alg:sspt}. Importantly, $\overline{\rho}$ is the lower-bound on the $\rho$ which leads to the satisfaction of the complementarity constraints. Our numerical experiments below reveal that this bound is very tight, and the selection in the order of $\rho=[1.2,2]\times\overline{\rho}$ makes the DCA solution strongly stationary for the original problem \eqref{prob:general}. 

\begin{algorithm}[t]
\caption{Selecting penalty parameter  $\rho$}
\label{alg:sspt}
\begin{algorithmic}[1]
    \REPEAT
        \STATE Sample $\tilde{\mb{d}}$ from the input domain $\mathcal{D}$
        \STATE Pass $\tilde{\mb{d}}$ through the trained NN to retrieve $\tilde{\mb{y}}$ and $\tilde{\mb{v}}$ 
    \UNTIL{$\tilde{y}_i \neq \tilde{v}_i$ for each neuron $i$}
    
    \STATE \textbf{return} $(\tilde{\mb{d}},\tilde{\mb{y}},\tilde{\mb{v}})$
    
    \STATE Obtain $(\overline{\mb{d}},\overline{\mb{y}},\overline{\mb{v}})$  from relaxation \eqref{prob:general_relax} on $(\tilde{\mb{d}},\tilde{\mb{y}},\tilde{\mb{v}})$
    
    \STATE Obtain duals $(\overline{\mbg{\lambda}}, \overline{\mbg{\mu}}^v,\overline{\mbg{\mu}}^y)$ by solving KKTs \eqref{KKTs_relax} as a feasibility problem with primal variables fixed to $(\overline{\mb{d}},\overline{\mb{y}},\overline{\mb{v}})$
    
    \STATE Compute the lower bound $\overline{\rho}$ from \eqref{eq:rho_bar}. 
    
    \STATE Fine-tune $\rho\geqslant\overline{\rho}$ from here until the DCA returns the solution with complementary conditions satisfied. 
\end{algorithmic}
\end{algorithm}

\section{Application to Data Centers Scheduling}\label{sec:application}

Our motivating problem is the optimal scheduling of electricity demand from large data centers in power systems. In this section, we first describe a popular yet difficult-to-implement approach to solving this problem as a bilevel optimization problem. We then detail an arguably more practical NN-based approach.

\subsection{Conventional Bilevel Formulation}

Consider an operator of spatially distributed data centers who allocates the workloads to minimize the cost of electricity consumption. This problem is solved in advance, when the operator receives a forecast of the total electricity demand, denoted by $\Delta$, to be allocated among grid nodes hosting data centers. The optimal allocation is driven by locational marginal prices (LMPs), defining the cost of electricity at specific nodes. The prices are computed via dual variables of the optimal power flow (OPF) problem \cite{chatzivasileiadis2018lecture}:
\begin{subequations}\label{problem:opf}
\begin{align}
\minimize{\underline{\mb{p}}\leqslant\mb{p}\leqslant\overline{\mb{p}}}\quad& \mb{p}^{\top}\mb{C}\mb{p}+\mb{c}^{\top}\mb{p} \label{opf_obj}\\[-1ex]
        \st\quad&\mb{1}^{\top}(\mb{p}-\mbg{\ell} - \mb{d})=0\quad\colon\lambda \label{opf_pb}\\
        &|\mb{F}(\mb{p}-\mbg{\ell} - \mb{d})|\leqslant\overline{\mb{f}}\quad\colon\underline{\mbg{\mu}},\overline{\mbg{\mu}}\label{opf_flow}
\end{align}
\end{subequations}
which seeks generator dispatch $\mb{p}$ that minimizes quadratic cost \eqref{opf_obj}. The dispatch must remain within a feasible range $[\underline{\mb{p}},\overline{\mb{p}}]$ and satisfy a set of grid constraints. Equality constraint \eqref{opf_pb} enforces a system-wide power balance between dispatched generation and loads composed of conventional loads $\mbg{\ell}$ and data center loads $\mb{d}$. The power flows in transmission lines are computed using the matrix of power transfer distribution factors  $\mb{F}$, which maps net power injections $(\mb{p}-\mbg{\ell}-\mb{d})$ to power flows as $\mb{F}(\mb{p}-\mbg{\ell}-\mb{d})$. They are capped by line capacity $\overline{\mb{f}}$ via \eqref{opf_flow}. The optimal dual variables (stated after the colon signs) are used to compute LMPs as
\begin{align}
    \mbg{\pi}=\mb{1}\lambda^{\star} - \mb{F}^{\top}\overline{\mbg{\mu}}^{\star} + \mb{F}^{\top}\underline{\mbg{\mu}}^{\star}
    \label{eq:LMP}
\end{align}
where the uniform price $\lambda^{\star}$ is corrected by line congestion terms, justifying the locational nature of electricity prices. 

To map a specific demand allocation $\mb{d}$ to LMPs, we need to solve the OPF sub-problem \eqref{problem:opf}. The optimal allocation becomes a bilevel optimization problem, which optimizes $\mb{d}$ using feedback from the embedded OPF sub-problem, i.e.,
\begin{subequations}\label{problem:bilevel}
\begin{align}
    \minimize{\mb{d}}\quad&\mbg{\pi}^{\top}\mb{d}\label{bilevel_obj}\\
    \st\quad&\underline{\mb{d}}\leqslant\mb{d}\leqslant\overline{\mb{d}} \label{bilevel_con1}\\
    &\mb{1}^{\top}\mb{d}=\Delta\label{bilevel_con2}\\
    &\mbg{\pi} \in \text{dual sol. of \eqref{problem:opf} formulated on $\mb{d}$} 
\end{align}
\end{subequations}
where the objective function \eqref{bilevel_obj}  minimizes the total charge for electricity consumption of data centers. Constraints in \eqref{bilevel_con1} ensure that the allocation respects the minimum and maximum capacity of data centers. The equality constraint \eqref{bilevel_con2} enforces the conservation of electricity demand. 


\subsection{Solution Based on the Optimization over Trained NN}

There are several practical challenges that prevent solving problem \eqref{problem:bilevel}. The first is the lack of data necessary to formulate OPF problem \eqref{problem:opf}, including actual network matrix $\mb{F}$, cost coefficients $\mb{C}$ and $\mb{c}$, and dispatch limits $\underline{\mb{p}}$ and $\overline{\mb{p}}$.  Second, solving \eqref{problem:bilevel} means optimizing over its KKTs that itself is an NP-hard problem \cite{dvorkin2018consensus}. Finally, the actual market engine can be more complicated than a well-posed convex formulation \eqref{problem:bilevel}, e.g., include unit commitment binary variables preventing its standard KKT reformulation. 

The optimization over trained NNs overcomes these challenges. NNs can approximate the mapping from data center loads $\mb{d}$ to resulting charges from historical data. Consider a labeled dataset  $\{(\mb{d}_{1},\lambda_{1}),\dots,(\mb{d}_{n},\lambda_{n})\}$, where for each load instance $\mb{d}_{i}$ we have the corresponding charge $\lambda_{i}$, obtained as the inner product of LMPs and data center loads from past outcomes. A trained NN, with fixed weight and bias parameters, takes input $\mb{d}$ and outputs $\lambda$. The optimization problem with the embedded trained NN takes the form:
\begin{subequations}\label{eq:dca_dc}
\begin{align}
   \minimize{\mb{d},\mb{y},\mb{v}, \lambda}
   \;\;&  \lambda\\
    \st\;\;
    &\underline{\mb{d}}\leqslant\mb{d}\leqslant\overline{\mb{d}} \\
    &\mb{1}^{\top}\mb{d}=\Delta\\
    & \mb{y}^1=\mb{W}^1 \mb{d}+\mb{b}^1 + \mb{v}^1,\\
    & \mb{y}^i=\mb{W}^i \mb{y}^{i-1}+\mb{b}^i + \mb{v}^i, \forall i=2,\dots,N \\
    &\mb{0}\leqslant \mb{y}^{i} \perp \mb{v}^{i} \geqslant \mb{0},\quad\quad\quad\!\;\forall i=1,\dots,N \\
     & \mbg{\lambda}= \mb{W}^{N + 1} \mb{y}^{N}+\mb{b}^{N+1},
    \end{align}
\end{subequations}
where the NN input takes the forecast $\Delta$ and computes the least-cost allocation of electricity demands $\mb{d}$. 

\section{Numerical Experiments}\label{sec:experiments}
We test the DCA algorithm on small and large networks. We relegate the modeling details to the online repository, including data and code to replicate the results:
\begin{center}\small \texttt{\href{https://github.com/xl359/NN_DC_OPT}{https://github.com/xl359/NN\_DC\_OPT}}.
\end{center}

\subsection{Experiments with Ground Truth Solution}
We first experiment with a small PJM $5$-bus system, where all 3 loads are assumed to be data centers with demand limits within the range of $[0.8,1]$  of the nominal load value. We train a NN with $2$ hidden layers, with $50$ neurons each, to compute the total demand charge on $10,000$ demand samples. For this small NN, we retrieve the global optimal solution to problem \eqref{eq:dca_dc} using \textit{Gurobi} solver. 

We apply the DCA Alg. \ref{alg:bilinear} to compute the optimal demand allocation, when the forecast $\Delta$ is $90$\% of the maximum total data center load. We used Alg. 2 to compute $\overline{\rho}=110.5$, and selected $\rho^{\star}=1.5\times\overline{\rho}$, which ensured the satisfaction of complementarity constraints; we use $\rho^\star$ as a reference.  

The ground truth solution yields the total demand charge of $\$194.4k$, which is matched by the DCA with $\rho^{\star}$ (with a negligible error of $<10^{-6})$. The trajectories of DCA solutions are given in Fig. \ref{fig:smallcasesmall}. Increasing the penalty leads to an order-of-magnitude increase in the number of iterations. Moreover, selecting a significantly larger penalty leads to the feasible yet sub-optimal solution, e.g., see $20\times\rho^{\star}$. The penalty parameters below the $\rho^{\star}$ threshold do not allow for the convergence to the ground truth. In fact, such trajectories, if converged, lead to points that do not satisfy the complementarity constraints and hence are of no use. The results confirm the convergence of Alg. \ref{alg:bilinear} to the ground truth on a small example and necessitate the application of Alg. \ref{alg:sspt} to properly tune the DCA penalty parameter. 
\begin{figure}
    \centering
    \includegraphics[width=0.48\textwidth]{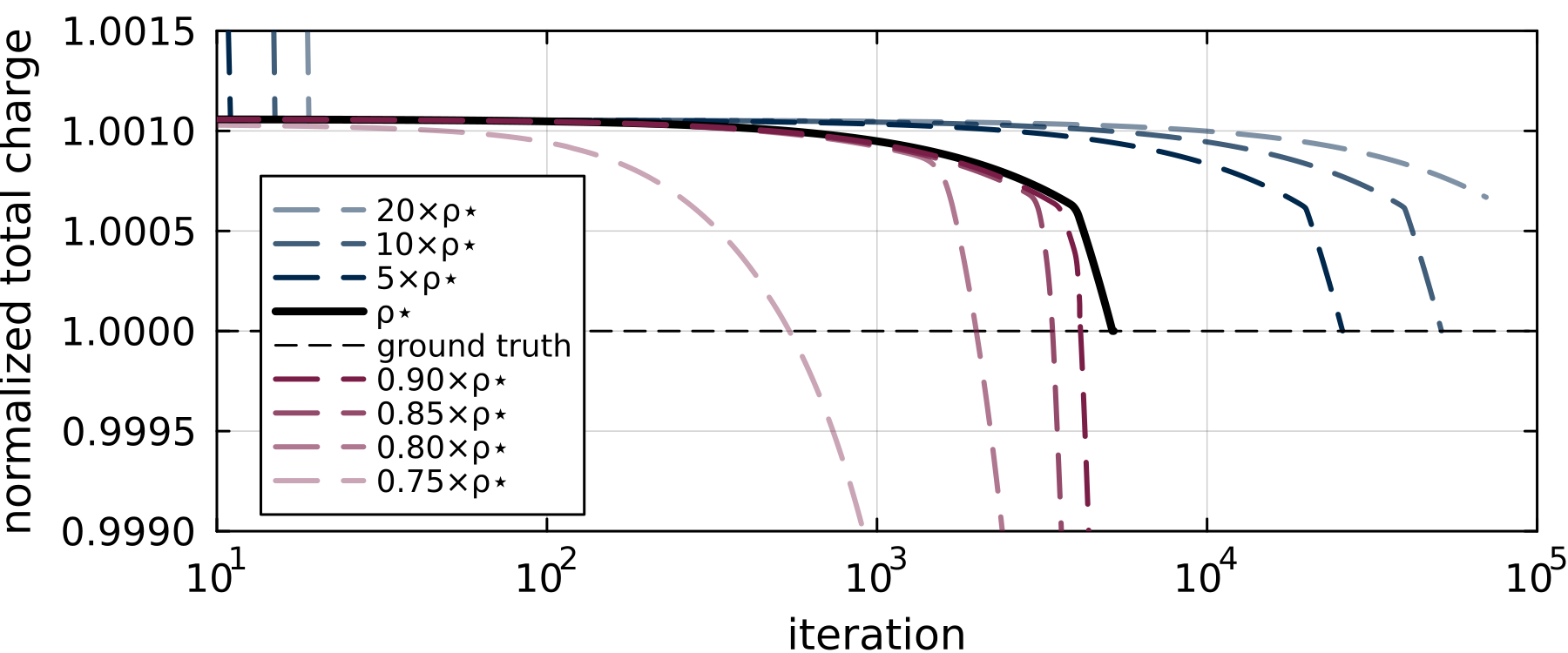} 
    \caption{5-bus PJM test case: DCA trajectories for varying penalty $\rho$. The total demand charge (cost of electricity) is normalized to the ground truth.}
    \label{fig:smallcasesmall}
\end{figure}

\subsection{Experiments on the IEEE 118-bus system}
We now consider the IEEE 118-bus system with 10 data centers installed, as shown in Fig. \ref{fig:network_lmp_graph}. The history of operations is modeled by $12,500$ samples of individual demands from the range of $[0,100]$ MWh, then map to the total demand charge. We approximate this mapping by a $5$-hidden layer NN with $1,000$ ReLU hidden neurons. We do not have the ground truth solution, as solving \eqref{eq:dca_dc} as a mixed-integer problem is computationally infeasible in this case. Instead, we benchmark against the baseline demand allocation sampled from the same distribution as the training samples. The DCA computes a new demand allocation in expectation of improving on the baseline. We randomly picked 50 baseline cases, of which DCA could successfully converge and roughly match the output of the OPF model on 47. We present the results for those cases. 

\begin{figure}
    \centering
    \includegraphics[width=0.48\textwidth]{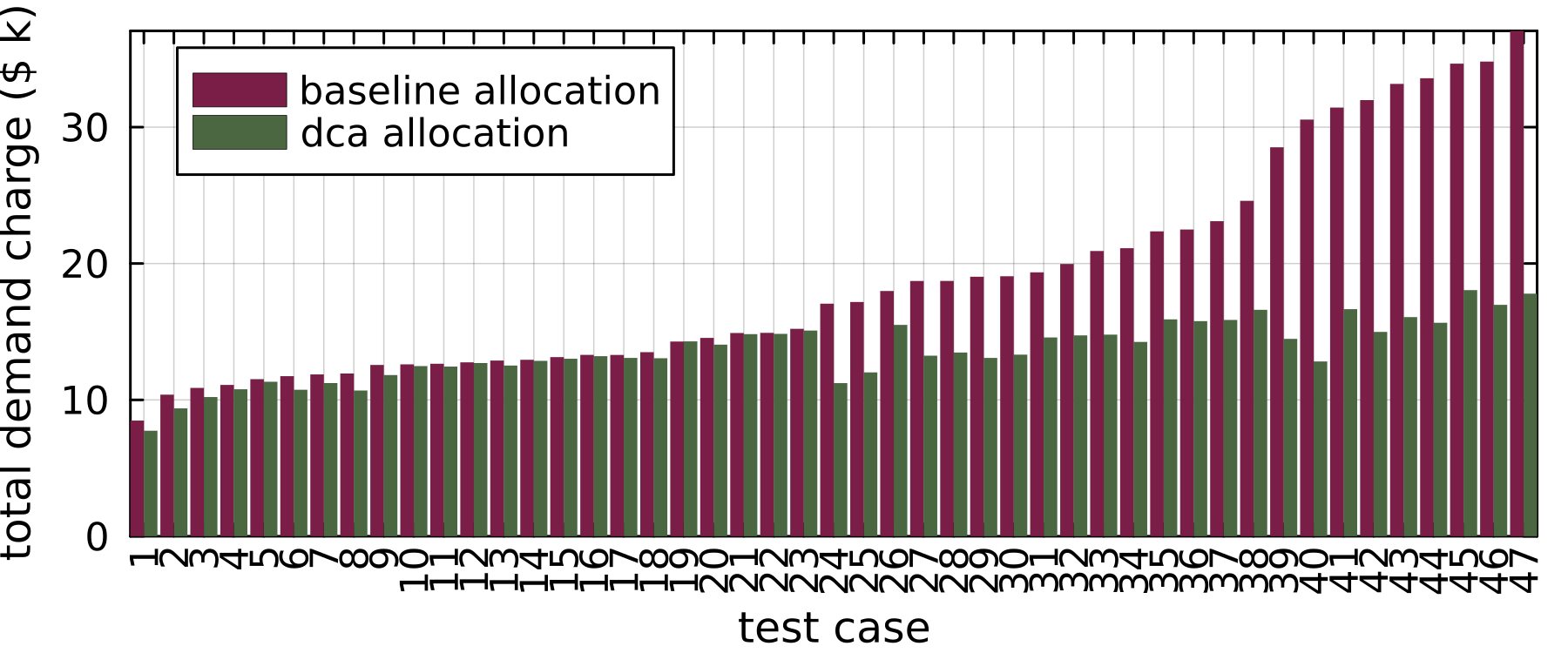}  
    \caption{Total charges for data center electricity consumption on the baseline and DCA allocations across 47 simulations in IEEE 118-bus system.}
    \label{fig:savemoney}
\end{figure}

\begin{figure}
    \centering
    \includegraphics[width=0.48\textwidth]{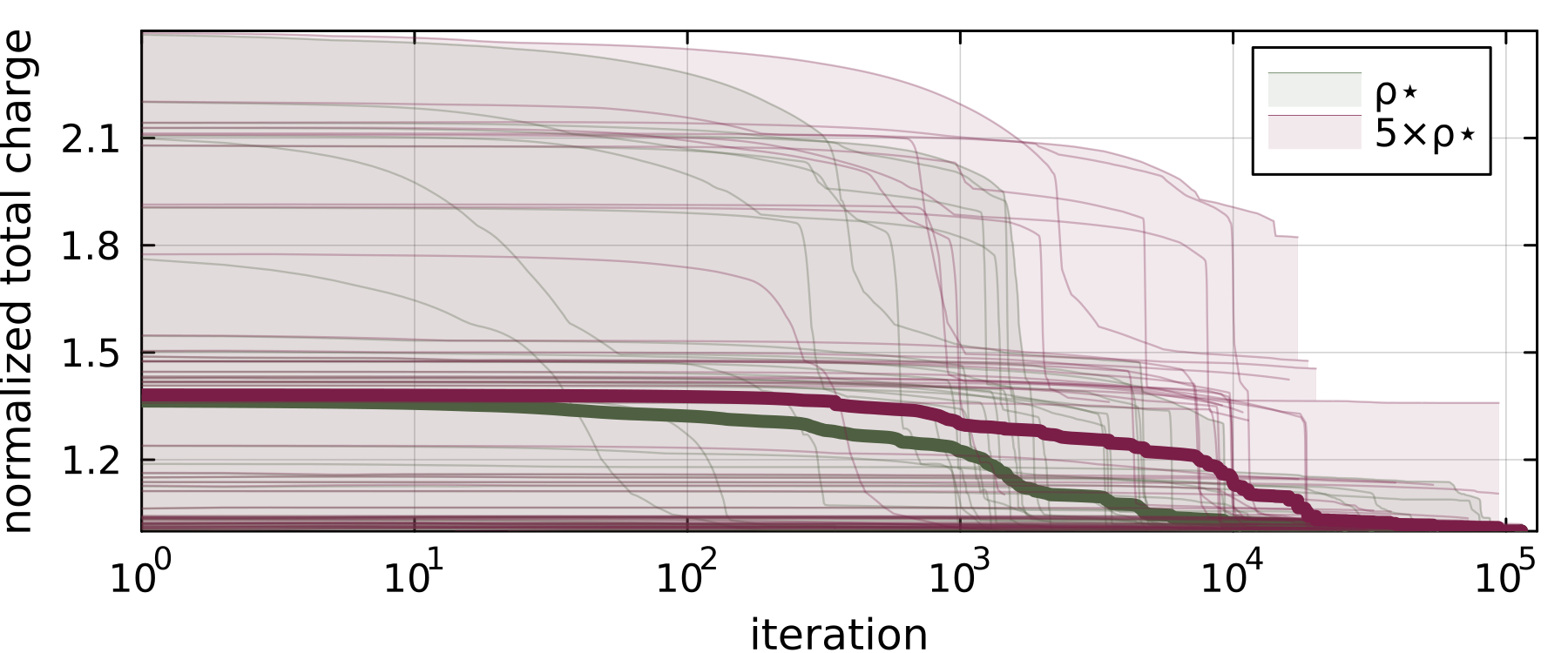}  
    \caption{DCA convergence in the experiments with the IEEE 118-bus system. The green thin lines indicates the actual trajectories using nominal penalty $\rho^\star$. The red thin lines indicates the convergence trajectories using $5\times\rho^\star$ penalty. The thick lines are the average trajectories. }
    \label{fig:convplot}
\end{figure}

\begin{figure}
    \centering
    \includegraphics[width=0.485\textwidth]{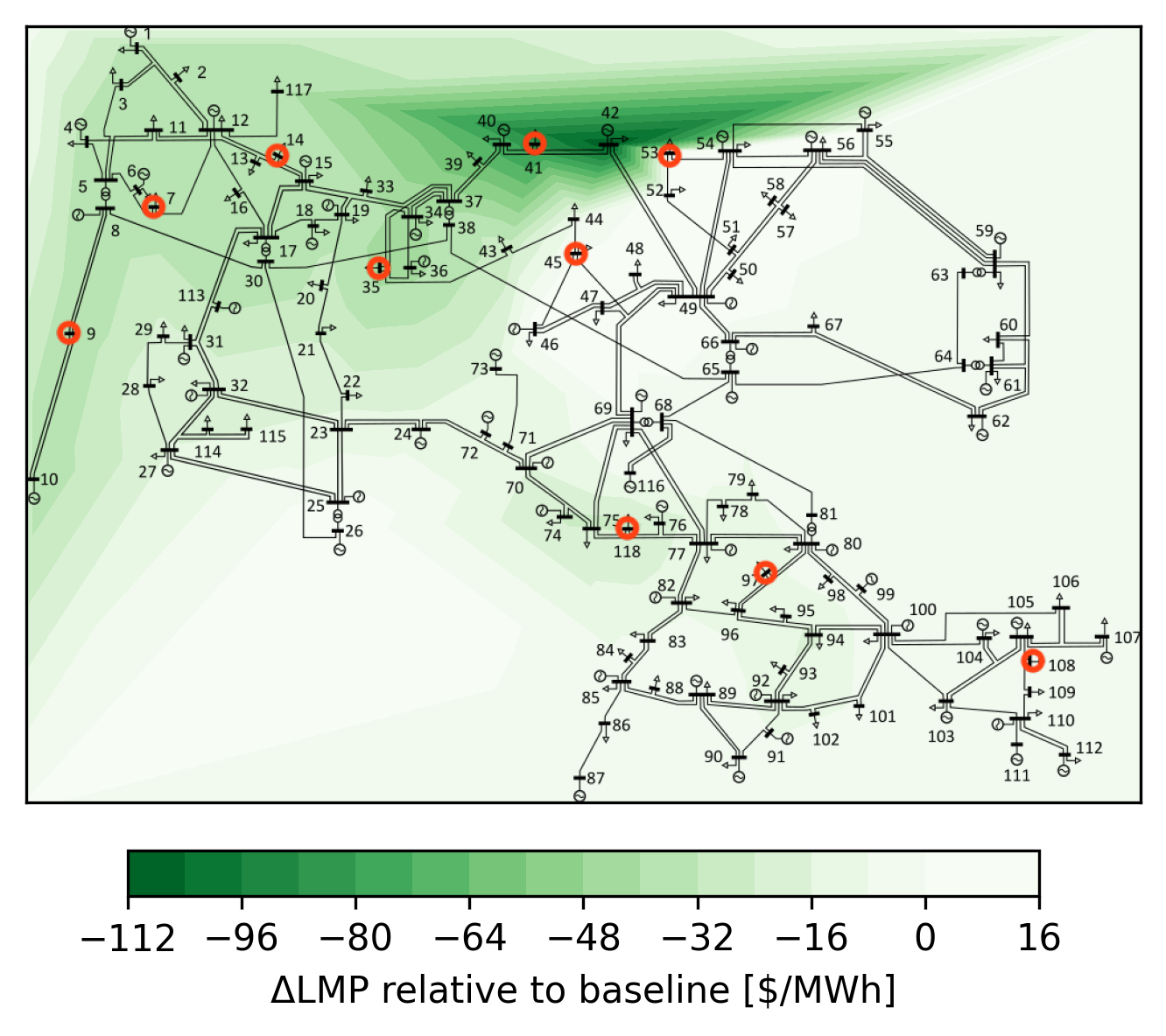}  
    \caption{The difference in LMPs induced by the baseline and DCA demand allocations in the IEEE 118-bus system. The red marks indicate buses hosting data centers. For the majority of these buses, the LMPs are significantly reduced under the DCA demand allocation compared to the baseline. Case $\#28$.}
    \label{fig:network_lmp_graph}
\end{figure}

Figure \ref{fig:savemoney} contrasts the charges resulted form the baseline and DCA demand allocations in $47$ cases, sorted by magnitude of the baseline allocation charge. In the first $23$ cases, the savings are small. These cases are lightly congested with almost uniform LMPs, where there is little to no benefit from the spatial load redistribution. In other cases, the DCA demand allocation significantly improves on the baseline. The DCA finds allocations that significantly reduces grid congestion. For one of such cases,  Fig. \ref{fig:network_lmp_graph} illustrates the difference between the LMPs induced on the baseline and the DCA solutions. At the majority of buses hosting data centers, the LMP is significantly reduced, e.g., up to $-85.4$ \$/MWh at bus $41$, leading to a significantly lower charges. 

Figure \ref{fig:convplot} presents the convergence behavior of DCA in $47$ cases. All values are normalized by the final DCA solution for the nominal $\rho^\star$ value. The smaller $\rho$ consistently leads to faster convergence. There is only one case out of the $47$ test cases for which the larger $\rho$ converged faster. The majority of cases converge within $10^5$ iterations, and on overage it takes below $10^4$ to converge. Each iteration takes $\approx 0.37$ seconds, allowing for convergence in around $1$ hour on average. 

Finally, we discuss some challenges in DCA convergence. There are a few cases where the DCA solution was different from that in the OPF solution on the same demand allocation. Such cases can be addressed by early stopping or keeping the best OPF solution so far in the DCA iterations. They can also be solved by increasing the penalty and stepping back several iterations. Only $3$ of $50$ cases were difficult to address and require further investigation.

\section{Conclusions}\label{sec:conclusion}

We develop a difference-of-convex algorithm to solve decision-making problems with objectives or constraints represented by trained neural networks. The algorithm avoids the computational bottleneck of ReLU logical constraints through informed penalization of such constraints in the objective function, hence avoiding computationally expensive mixed-integer techniques. The application to the problem of the optimal data center demand allocation in power grids revealed a large reduction in the total demand charge relative to the baseline demand allocation. Our results justify applications in the settings with a clear baseline solution upon which the algorithm can improve. As future research, we will expand the scope of modeling constraints, e.g., to include both spatial and temporal demand allocation, and provide more theoretical guarantees on the convergence behavior.

\bibliographystyle{IEEEtran}
\bibliography{IEEEabrv,reference}

\endgroup
\end{document}